\documentclass[11pt]{article}
\usepackage{latexsym,amsmath,amsfonts,url}
\usepackage{graphicx,subfigure}
\newcommand{\ignore}[1]{}

\begin{document}
\begin{center}
{\LARGE\textbf{Bipancyclic subgraphs in random bipartite graphs}}\\
\bigskip

Yilun Shang\\
Institute for Cyber Security, University of Texas at San Antonio\\San Antonio, Texas 78249, USA\\
Email: \texttt{shylmath@hotmail.com}

\end{center}

\begin{abstract}

A bipartite graph on $2n$ vertices is bipancyclic if it contains
cycles of all even lengths from $4$ to $2n$. In this paper we prove
that the random bipartite graph $G(n,n,p)$ with $p(n)\gg n^{-2/3}$
asymptotically almost surely has the following resilience property:
Every Hamiltonian subgraph $G'$ of $G(n,n,p)$ with more than
$(1/2+o(1))n^2p$ edges is bipancyclic. This result is tight in two
ways. First, the range of $p$ is essentially best possible. Second,
the proportion $1/2$ of edges cannot be reduced. Our result extends
a classical theorem of Mitchem and Schmeichel.

\bigskip

\end{abstract}

\bigskip

\section{Introduction}

A bipartite graph on $2n$ vertices is called bipancyclic if it
contains cycles of all even lengths from $4$ to $2n$. Analogously, a
graph on $n$ vertices is called pancyclic if it contains cycles of
all length $t$ for $3\le t\le n$. Clearly, (bi)pancyclic graphs are
Hamiltonian but the converse is not true in general. A variety of
sufficient conditions for a Hamiltonian bipartite graph to be
bipancyclic have been studied in the literature, including
\cite{5,3,6,4} and \cite{2}. Recall that a bipartite graph is called
balanced if the two classes of bipartition have the same
cardinality. In \cite{6} Mitchem and Schmeichel proved the following
theorem.

\smallskip
\noindent\textbf{Theorem 1.} \quad \itshape Let $G$ be a Hamiltonian
bipartite balanced graph with $2n$ vertices and $m$ edges. If
$m>n^2/2$, then $G$ is bipancyclic. \normalfont
\smallskip

Recently, Sudakov and Vu \cite{7} proposed the framework of
resilience of graphs, in which many extremal graph-theoretic
properties such as Hamiltonicity and pancyclicity can be studied
(see e.g. \cite{9,10,8,1,11}). Let $\mathcal{P}$ be a monotone
increasing graph property. Define the global resilience of a graph
$G$ with respect to $\mathcal{P}$ as the minimum number $r$ such
that by deleting $r$ edges from $G$, one can obtain a graph not
having $\mathcal{P}$. Using this notion, the above Theorem 1 can be
reformulated as a global resilient statement with an additional
constraint: If one deletes fewer than $n^2/2$ edges from the
complete bipartite graph $K_{n,n}$ while preserving Hamiltonicity,
then the resulting graph is always bipancyclic.

In this paper, we study bipancyclicity of random bipartite graphs in
the context of global resilience by extending Theorem 1. The model
of random bipartite graphs $G(n,n,p)$ is the probability
distribution on the set of all bipartite balanced graphs with vertex
set $\{1,2,\cdots,2n\}$ such that each pair of vertices from
different classes of bipartition forms an edge randomly and
independently with probability $p$. Its monopartite version is the
celebrated binomial random graph $G(n,p)$ (see e.g. \cite{12}). We
say that $G(n,n,p)$ (or $G(n,p)$) possesses a graph property
$\mathcal{P}$ asymptotically almost surely, or a.a.s. for short, if
the probability that $G(n,n,p)$ (or $G(n,p)$) possesses
$\mathcal{P}$ approaches to 1 as $n$ tends to infinity. Lee and
Samotij \cite{1} recently proved that if $p\gg n^{-1/2}$, then
$G(n,p)$ a.a.s. satisfies the following: Every Hamiltonian subgraph
$G'\subset G(n,p)$ with more than $(\frac12+o(1))n^2p/2$ edges is
pancyclic. Our main result is a corresponding version for
$G(n,n,p)$, which is the following generalization of Theorem 1
(since $G(n,n,1)=K_{n,n}$).

\smallskip
\noindent\textbf{Theorem 2.} \quad \itshape If $p\gg n^{-2/3}$, then
$G(n,n,p)$ a.a.s. satisfies the following. Every Hamiltonian
subgraph $G'\subset G(n,n,p)$ with more than $(1+o(1))n^2p/2$ edges
is bipancyclic. \normalfont
\smallskip

Theorem 2 is asymptotically tight in two ways. First, one cannot
improve the exponent $-2/3$. To see this, assume that $p\ll
n^{-2/3}$ and fix a Hamilton cycle $H$ in $G(n,n,p)$. From each
4-cycle in $G(n,n,p)$, delete one edge which does not belong to $H$.
Since a.a.s there are at most $n^4p^4=o(n^2p)$ 4-cycles in the
graph, only a small proportion of edges is deleted and the resulting
graph does not contain any 4-cycles, hence not bipancyclic. Second,
Hamilton subgraphs with fewer than $(1+o(1))n^2p/2$ edges need not
be bipancyclic. Assume that $p\gg n^{-2/3}$ and fix a Hamilton cycle
$H$ in $G(n,n,p)$. We label the vertices as shown in Fig. 1 such
that $H=\{0,1,2,\cdots,2n-1,0\}$. Delete all edges $\{0,j\}$ from
$G(n,n,p)$ except two edges $\{0,1\}$ and $\{0,2n-1\}$. For each
even $i$ with $2\le i\le 2n-2$, delete all edges $\{i,j\}$ from
$G(n,n,p)$ with $j\ge i+3$. A.a.s. we will delete at most
$(1+o(1))n^2p/2$ edges, and a.a.s. the above process produces a
graph $G'$ with at least $(1+o(1))n^2p/2$ edges. Note that $H\subset
G'$ so $G'$ is Hamiltonian. However, $G'$ contains no 4-cycles, thus
not bipancyclic.

\begin{figure}[!t]
\begin{center}
\scalebox{0.7}{\includegraphics[166pt,317pt][447pt,557pt]{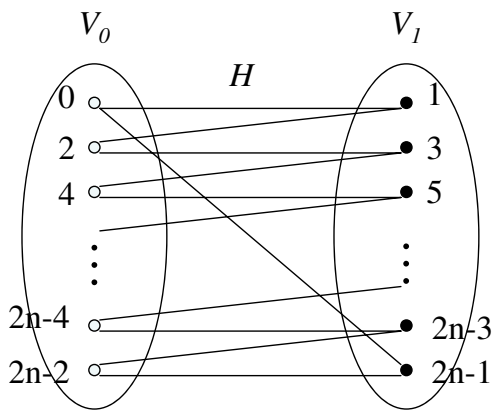}}\caption{A
bipartite balanced graph with bipartition $|V_0|=|V_1|=n$.
$H=\{0,1,2,\cdots,2n-1,0\}$ is a Hamilton cycle.}
\end{center}
\end{figure}

The rest of the paper is organized as follows. In Section 2 we
present some notations and preliminaries that will be needed in our
development later. The proof of Theorem 2 comprises two parts: In
Section 3 we establish the existence of short and long cycles of
even lengths, while in Section 4 we establish the existence of
medium ones of even lengths.

\section{Preliminaries}

Let $G=(V,E)$ denote a graph with vertex set $V$ and edge set $E$.
Similarly, a bipartite graph $G$ with edge set $E$ is denoted by
$G=(V_0,V_1,E)$ where $V_0$ and $V_1$ are the two classes of the
bipartition. For a vertex $v$, we denote its neighborhood by $N(v)$,
and its degree by $\operatorname{deg}(v)=|N(v)|$. For a set $X$, let
$E(X)$ be the set of edges in the induced subgraph $G[X]$, and let
$e(X)=|E(X)|$. When we have several graphs under consideration, we
may use subscripts such as $\operatorname{deg}_G(v)$ to indicate the
graph we are currently working with. We often omit floor and ceiling
signs whenever these are not crucial. We also assume the order of
the graphs is large enough throughout our derivation.

The following concentration inequality (see e.g. \cite[Corollary
2.3]{12}) will be often used in the proof of main result.

\smallskip
\noindent\textbf{Theorem 3.} (Chernoff's inequality) \quad \itshape
Let $0<\varepsilon\le3/2$. If $X$ is a binomial random variable with
parameter $n$ and $p$, then
$$
P(|X-\mathbb{E}(X)|\ge\varepsilon\mathbb{E}(X))\le2e^{-\varepsilon^2\mathbb{E}(X)/3},
$$
where $\mathbb{E}$ represents the expectation operator. \normalfont
\smallskip

The following results on cycles of fixed length were proved by Bondy
and Simonovits \cite{13}, and Haxell et al. \cite{13a}. They yield
the existence of very short cycles.

\smallskip
\noindent\textbf{Theorem 4.} \quad \itshape (i) \cite{13} Let $k$ be
a positive integer, and let $G$ be a graph on $n$ vertices with more
than $100kn^{1+(1/k)}$ edges. Then $G$ contains a cycle of length
$2k$.

(ii) \cite{13a} For any fixed integer $l\ge2$ and
$\varepsilon\in(0,1)$, there exists a constant $C>0$ such that if
$p\ge Cn^{-1+1/(2l-1)}$, then $G(n,p)$ a.a.s. satisfies the
following. Every subgraph $G'\subset G(n,p)$ with at least
$(1+\varepsilon)n^2p/8$ edges contains a cycle of length $2l$.
\normalfont
\smallskip

As reasoned in \cite{1}, our proof of Theorem 2 will rely on a
hypergraph construction, which fits in the general framework
developed for extremal properties of random discrete structures by
Schacht \cite{14} (similar results were obtained by Conlon and
Gowers \cite{15} independently). Before introducing the general
theorem, we need some definitions.

\smallskip
\noindent\textbf{Definition 1.} \quad \itshape Let $H$ be a
$k$-uniform hypergraph, $\alpha\ge0$ and $\varepsilon_0\in(0,1)$.
Let $f:(0,1)\rightarrow(0,1)$ be a non-decreasing function. We say
$H$ is $(\alpha,f,\varepsilon_0)$-dense if the following is true.

For every $\varepsilon\ge\varepsilon_0$ and every $U\subset V(H)$
with $|U|\ge(\alpha+\varepsilon)|V(H)|$, we have
$$
|E(H[U])|\ge f(\varepsilon)|E(H)|.
$$
\normalfont
\smallskip

For a $k$-uniform hypergraph $H$, $v\in V(H)$, $U\subset V(H)$ and
$i\in\{1,2,\cdots,k-1\}$, we define
$$
\operatorname{deg}_i(v,U)=|\{X\in E(H): |X\cap
(U\backslash\{v\})|\ge i\ \mathrm{and}\ v\in X\}|.
$$
For $q\in[0,1]$ and a set $X$, let $X_q$ be the binomial random
subset of $X$ with survival probability $q$.

\smallskip
\noindent\textbf{Definition 2.} \quad \itshape Let $H$ be a
$k$-uniform hypergraph, $p\in(0,1)$ and $K\ge1$. We say $H$ is
$(K,p)$-bounded if the following is true.

For every $i\in\{1,2,\cdots,k-1\}$ and $q\in[p,1]$, we have
$$
\mathbb{E}\left(\sum_{v\in
V(H)}\operatorname{deg}_i(v,V(H)_q)^2\right)\le
Kq^{2i}\frac{|E(H)|^2}{|V(H)|}.
$$
\normalfont
\smallskip

\smallskip
\noindent\textbf{Theorem 5.} (\cite{14}) \quad \itshape Suppose
$(H_n)_{n\in\mathbb{N}}$ is a sequence of $k$-uniform hypergraphs.
Let $(p_n)_{n\in\mathbb{N}}$ be a sequence of probabilities. Let
$(v_n)_{n\in\mathbb{N}}$ and $(e_n)_{n\in\mathbb{N}}$ be sequences
of integers satisfying $p_nv_n\rightarrow\infty$ and
$p_n^ke_n\rightarrow\infty$ as $n\rightarrow\infty$. Let
$\alpha\ge0$, $K\ge1$ and $f:(0,1)\rightarrow(0,1)$ be a
non-decreasing function. For every $\varepsilon\in(0,1)$, there
exist $\varepsilon_0\in(0,1)$, $b>0$, $C\ge1$ and $n_0\ge1$ such
that for every $n\ge n_0$ and every $q$ with $n^{-1/3}\ge q\ge Cp_n$
the following holds.

If $H_n$ is $(\alpha,f,\varepsilon_0)$-dense and $(K,p_n)$-bounded
satisfying $|V(H_n)|\ge v_n$ and $|E(H_n)|\ge e_n$, then with
probability at least $1-e^{-bqv_n}$, every subset $W\subset
V(H_n)_q$ with $|W|\ge(\alpha+\varepsilon)|V(H_n)_q|$ contains an
edge of $H_n$. \normalfont
\smallskip

\section{Existence of short and long cycles of even lengths}

In this section, we establish the existence of short and long cycles
of even lengths in the following sense.

\smallskip
\noindent\textbf{Theorem 6.} \quad \itshape For any
$\varepsilon\in(0,1)$, there exist constants $C>0$ and
$\delta\in(0,1)$ such that if $p\ge Cn^{-2/3}$, then $G(n,n,p)$
a.a.s. satisfies the following. Every Hamiltonian subgraph
$G'\subset G(n,n,p)$ with more than $(1+\varepsilon)n^2p/2$ edges
contains a cycle of length $t$ for all even $t\in[4,2\delta
n]\cup[2(1-\delta)n,2n]$. \normalfont
\smallskip

We will follow the idea of \cite{1} and separate the proof of
Theorem 2 into two parts: Theorem 6 is responsible for short and
long cycles, and Theorem 7 (see Section 4 below) will be responsible
for intermediate cycles. Note that if we choose $C$ in Theorem 6 to
be large enough, the existence of 4-cycle and 6-cycle follows easily
from Theorem 4. Hence, in what follows we will focus on cycles of
length 8 and above.

Let $[2n]$ be the set of remainders modulo $2n$, namely
$[2n]=\{0,1,2,\cdots,2n-1\}$. The addition of the elements of $[2n]$
will be performed modulo $2n$ throughout this paper. A labeling of
the vertices of the complete bipartite graph $K_{n,n}=(V_0,V_1,E)$
is called allowable if the vertices in $V_0$ are labeled by even
numbers and those in $V_1$ are labeled by odd ones; c.f. Fig. 1. Fix
an allowable labeling and let $C_{2n}$ be the subgraph of $K_{n,n}$
consisting of the edges $\{i,i+1\}$ for all $i\in[2n]$. For
illustration we may draw it as a circle with labels
$0,1,2,\cdots,2n-1$ in the clockwise order. For each $i\in[2n]$,
denote its distance from 0 on the cycle $C_{2n}$ by $\|i\|$, namely
$$
\|i\|=\min\{k\ge0: k\equiv i \mod 2n\ \mathrm{or}\ k\equiv -i \mod
2n\}.
$$
For an even $l$ with $0\le l\le n$, a 4-vertex subgraph $X\subset
K_{n,n}\backslash C_{2n}$ is called an $l$-shortcut if it is of one
of the following types:
\begin{itemize}
\item[(i)] There are $i_1,i_2,i_3,i_4\in[2n]$ such that
$i_1,i_1+1,i_2,i_2+1,i_3,i_3+1,i_4$ and $i_4+l+1$ are all distinct
and lie in the clockwise order on $C_{2n}$, and $X$ is composed of
the edges $\{i_1,i_3\}$, $\{i_1+1,i_4\}$, $\{i_2,i_4+l+1\}$ and
$\{i_2+1,i_3+1\}$. Moreover, $i_1+1$ and $i_2$ belong to different
classes of bipartition. So do $i_1$ and $i_2+1$.

\item[(ii)]  There are $i_1,i_2,i_3,i_4\in[2n]$ such that
$i_1,i_1+1,i_2,i_2+1,i_4,i_4+l+1,i_3$ and $i_3+1$ are all distinct
and lie in the clockwise order on $C_{2n}$, and $X$ is composed of
the edges $\{i_1,i_3\}$, $\{i_1+1,i_4\}$, $\{i_2,i_4+l+1\}$ and
$\{i_2+1,i_3+1\}$. Moreover, $i_1+1$ and $i_2$ belong to different
classes of bipartition. So do $i_1$ and $i_2+1$.
\end{itemize}

Since our formulation is similar with that in \cite{1}, we highlight
the novelties in our methodology.
\begin{itemize}
\item The introduction of allowable labeling make it possible to
extend the analysis from monopartite graph $G(n,p)$ to bipartite
graph $G(n,n,p)$.

\item Additional constraints are posed on the definition of
$l$-shortcut. These modifications imply that $C_{2n}$ and $X$ may be
concurrent under any allowable labeling and will be critical in the
proof of some technical lemmas later.
\end{itemize}

A key observation is that the graph $C_{2n}\cup X$ contains cycles
of lengths $l+8$ and $2n-l$ for every even $l$ with $0\le l\le n$
and every $l$-shortcut $X$; see Fig. 2. We will prove Theorem 6
through a series of lemmas.

\begin{figure}[!t]
\begin{center}
\scalebox{0.7}{\includegraphics[89pt,339pt][567pt,540pt]{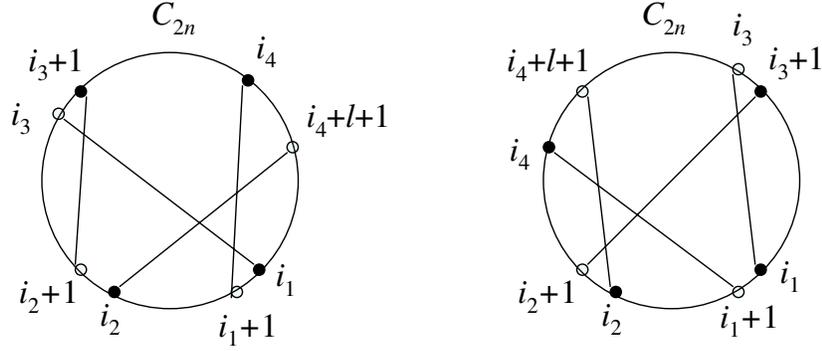}}\caption{Two
$l$-shortcuts of types (i) and (ii), respectively.}
\end{center}
\end{figure}

\smallskip
\noindent\textbf{Lemma 1.} \quad \itshape For any
$\varepsilon_0\in(0,1)$, there exist an $n_0\ge1$ such that if
$\varepsilon'\ge\varepsilon_0$ and $n\ge n_0$, then every
$2n$-vertex bipartite balanced graph $G'$ with
$e(G')\ge(1+\varepsilon')n^2/2$ contains at least
$\varepsilon'^8n^4/(4\cdot16^7)$ many of $l$-shortcuts for every
even $l$ with $0\le l\le \varepsilon'n/8$ and every allowable
labeling of the vertex set of $G'$ with $[2n]$. \normalfont
\smallskip

\noindent\textbf{Proof.} Assume that $\varepsilon'\ge\varepsilon_0$
and $n\ge n_0=192/\varepsilon_0$. Fix an allowable labeling of the
vertex set of $G'$ with $[2n]$. We have
$$
\sum_{i=0}^{n-1}\left(\deg_{G'}(2i)+\deg_{G'}(2i+1)\right)=2e(G')\ge(2+\varepsilon')\frac{n^2}{2}.
$$
Define a set $I=\{i\in\{0,1,\cdots
n-1\}:\deg_{G'}(2i)+\deg_{G'}(2i+1)\ge(1+\varepsilon'/2)n\}$. Via a
simple proof by contradiction we can see that
$|I|\ge\varepsilon'n/2$. For every $k$ with $0\le k\le
(1-\varepsilon'/4)n$, define $I(k)=\{i\in
I:\deg_{G'}(2i)\in[k,k+\varepsilon'n/4)\}$. Again applying proof by
contradiction we obtain that there exists some $k$ such that
$$
|I(k)|\ge\frac{|I|}{\lceil\frac{4}{\varepsilon'}\rceil}\ge\frac{\varepsilon'|I|}{8}\ge\frac{\varepsilon'^2}{16}n.
$$
We define $I'=I(k)$ for any such $k$. Therefore, for all $i,j\in I'$
we have
\begin{equation}
\deg_{G'}(2i)+\deg_{G'}(2j+1)\ge\deg_{G'}(2j)+\deg_{G'}(2j+1)-\frac{\varepsilon'}{4}n\ge\left(1+\frac{\varepsilon'}{4}\right)n.\label{1}
\end{equation}

On the other hand, it is easy to see that there exists a subset
$I''\subset I'$ satisfying $|I''|\ge\varepsilon'|I'|/32$, such that
for all $i,j\in I''$, the distance between $2i$ and $2j$ on the
cycle $C_{2n}$ satisfies $\|2i-2j\|\le\varepsilon'n/16$. We claim

\smallskip
\noindent\textbf{Claim 1.} \quad \itshape For all different $i,j\in
I''$ and every even $l$ with $0\le l\le \varepsilon'n/8$, there are
at least $(\varepsilon'n/32)^2$ many of $l$-shortcuts with
$\{i_1,i_2\}=\{2i,2j\}$. \\ \normalfont If this is true, the total
number of $l$-shortcut in $G'$ is at least
$$
{|I''|\choose2}\left(\frac{\varepsilon'n}{32}\right)^2\ge\frac{|I''|^2\varepsilon'^2}{4\cdot32^2}n^2\ge\frac{|I'|^2\varepsilon'^4}{4\cdot32^4}n^2
\ge\frac{\varepsilon'^8}{4\cdot 16^7}n^4.
$$

What remains to prove is Claim 1. Fix an even $l$ with $0\le l\le
\varepsilon'n/8$. Note that $\|2i-2j\|\le\varepsilon'n/16$.
Therefore, we may set $\{i_1,i_2\}=\{2i,2j\}$ such that $i_2=i_1+k$
for some $0<k\le\varepsilon'n/8$. Let $A$ be the set
$\{i_2+2,\cdots,i_1-1\}$ of vertices of $C_{2n}$ lying on the major
arc connecting $i_2+1$ to $i_1$. Let $A'=\{i\in A: i+l+1\in A\}$.
Therefore, we have
$$
|A'|=|A|-(l+1)\ge\left(2-\frac{\varepsilon'}{16}\right)n-2-l-1\ge\left(2-\frac{3\varepsilon'}{16}-\frac{\varepsilon'}{32}\right)n.
$$
Let $B=\{i\in A':\{i_1+1,i\},\{i_2,i+l+1\}\in E(G')\}$,
$N_1=N_{G'}(i_1+1)$ and $N_2=\{i\in[2n]:i+l+1\in N_{G'}(i_2)\}$.
Employing (\ref{1}) and the fact that $B=A'\cap N_1\cap N_2$, we
obtain
\begin{eqnarray*}
\left(1+\frac{\varepsilon'}{4}\right)n&\le&\deg_{G'}(i_1+1)+\deg_{G'}(i_2)\\
&=&|N_1\cup N_2|+|N_1\cap N_2|\\
&=&|N_1\cup N_2|+|([2n]\backslash A')\cap N_1\cap N_2|+|A'\cap
N_1\cap N_2|\\
&\le&n+2n-|A'|+|B|,
\end{eqnarray*}
where the last inequality holds since $l$ is even and $i_1+1$ and
$i_2$ belong to different classes of bipartition. Hence
\begin{equation}
|B|\ge\left(\frac{\varepsilon'}{4}-\frac{3\varepsilon'}{16}-\frac{\varepsilon'}{32}\right)n=\frac{\varepsilon'}{32}n.\label{2}
\end{equation}
Fix some $i_4\in B$, and let $J=\{i_4,\cdots,i_4+l+1\}$, $A''=\{i\in
A\backslash J:i+1\in A\backslash J\}$. Therefore, we have
$$
|A''|=|A|-|J|-2\ge\left(2-\frac{3\varepsilon'}{16}-\frac{\varepsilon'}{32}\right)n.
$$
Let $D=\{i\in A'':\{i_1,i\},\{i_2+1,i+1\}\in E(G')\}$. We can argue
analogously as above to derive
\begin{equation}
|D|\ge\left(\frac{\varepsilon'}{4}-\frac{3\varepsilon'}{16}-\frac{\varepsilon'}{32}\right)n=\frac{\varepsilon'}{32}n.\label{3}
\end{equation}
Fix some $i_3\in D$ and we readily have an $l$-shortcut $X\subset
G'$ consisting of edges
$\{i_1,i_3\},\{i_1+1,i_4\},\{i_2,i_4+l+1\},\{i_2+1,i_3+1\}$.
Therefore, the claim follows from (\ref{2}) and (\ref{3}). $\Box$

A key ingredient of the proof of Theorem 6 is to construct a
hypergraph to which Theorem 5 can be applied. Inspired by the
construction presented in \cite{1}, we define $H_n^l$ be a 4-uniform
hypergraph with the vertex set $V(H_n^l)=E(K_{n,n})$ for every
integer $n$ and every even $l$ with $0\le l\le n$. The hyperedges of
$H_n^l$ consist of all $l$-shortcuts in $K_{n,n}\backslash C_{2n}$.
Therefore, we have
\begin{equation}
|V(H_n^l)|=|E(K_{n,n})|=n^2\quad \mathrm{and}\quad
cn^4\le|E(H_n^l)|\le2n^4,\label{4}
\end{equation}
where $c>0$ is an absolute constant. The following corollary is
immediate from Lemma 1 and Definition 1.

\smallskip
\noindent\textbf{Corollary 1.} \quad \itshape Let
$f:(0,1)\rightarrow(0,1)$ be the function defined by
$f(\varepsilon')=4\varepsilon'^8/16^6$ for all
$\varepsilon'\in(0,1)$. For any $\varepsilon_0\in(0,1)$, there exist
$\delta\in(0,1)$ and $n_1\ge1$ such that for all $n\ge n_1$ and even
$l$ with $0\le l\le 2\delta n$, the hypergraph $H_n^l$ is
$(1/2,2f,\varepsilon_0/2)$-dense. \normalfont
\smallskip

Using an almost identical argument of \cite[Lemma 3.6]{1} we can
establish the $(K,n^{-2/3})$-boundedness of $H_n^l$. We leave the
proof of the following result to the reader.

\smallskip
\noindent\textbf{Lemma 2.} \quad \itshape There exists a constant
$K>0$ such that for all integer $n$ and even $l$ with $0\le l\le n$,
the hypergraph $H_n^l$ is $(K,n^{-2/3})$-bounded. \normalfont
\smallskip

For $\delta\in(0,1)$, define a monotone increasing graph property
$\mathcal{P}_{\delta}$ as follows. A $2n$-vertex bipartite balanced
graph $G$ satisfies $\mathcal{P}_{\delta}$ if and only if $G$
contains an $l$-shortcut for every even $l$ with $0\le l\le 2\delta
n$ and every allowable labeling of the vertices of $G$ with $[2n]$.

\smallskip
\noindent\textbf{Lemma 3.} \quad \itshape For any
$\varepsilon\in(0,1)$, there exist $\delta\in(0,1)$ and $C>0$ such
that if $Cn^{-2/3}\le p\le n^{-1/3}$, then $G(n,n,p)$ a.a.s.
satisfies the following. Every subgraph $G'\subset G(n,n,p)$ with
more than $(1/2+\varepsilon/2)e(G(n,n,p))$ edges satisfies
$\mathcal{P}_{\delta}$. \normalfont
\smallskip

\noindent\textbf{Proof.} Set $k=4$, $\alpha=1/2$, $p_n=n^{-2/3}$,
$v_n=n^2$ and $e_n=cn^4$ with $c$ given in (\ref{4}). Let $f$ be the
function defined in Corollary 1 and $K$ be given in Lemma 2. We have
$e_n(n^{-2/3})^4\rightarrow\infty$ and
$v_nn^{-2/3}\rightarrow\infty$, as $n\rightarrow\infty$.
Furthermore, let $\varepsilon_0,b,C,n_0$ be the numbers satisfying
the conclusion of Theorem 5. Let $\delta$ and $n_1$ be the numbers
given in Corollary 1 by using the parameter $\varepsilon_0$.

Suppose that $n\ge\max\{n_0,n_1\}$ and fix an even $l$ with $0\le
l\le 2\delta n$. Fix an allowable labeling of the vertices of
$G(n,n,p)$ with $[2n]$. In view of Corollary 1 and Lemma 2 we
observe that $H_n^l$ is $(1/2,2f,\varepsilon_0/2)$-dense and
$(K,n^{-2/3})$-bounded. Let $Cn^{-2/3}\le p\le n^{-1/3}$. Since
$|V(H_n^l)|=v_n$, $|E(H_n^l)|\ge e_n$ and $n\ge n_0$, Theorem 5
implies that with probability at least $1-e^{-bpn^2}$, every
subgraph $G'\subset G(n,n,p)$ with $e(G')\ge
(1/2+\varepsilon/2)e(G(n,n,p))$ contains a hyperedge of $H_n^l$,
which is an $l$-shortcut with respect to the above fixed labeling.
An application of Stirling's approximation shows that with
probability at least $1-(n!)^2ne^{-bpn^2}=1-o(1)$, the random
bipartite graph $G(n,n,p)$ satisfies the conclusion of Lemma 3.
$\Box$

\noindent\textbf{Completion of the proof of Theorem 6.} Let $\delta$
and $C$ be the numbers satisfying the conclusion of Lemma 3 by using
the parameter $\varepsilon/4$ (instead of $\varepsilon$). Let
$p'=Cn^{-2/3}$. An application of Theorem 3 shows that
$e(G(n,n,p'))\le(1+\varepsilon/8)n^2p'$ a.a.s. Hence, by using Lemma
3 we obtain that a.a.s. every subgraph of $G(n,n,p')$ with more than
$(1/2+\varepsilon/4)n^2p'$ edges satisfies $\mathcal{P}_{\delta}$.
We claim

\smallskip
\noindent\textbf{Claim 2.} \quad \itshape Assume that $0<p'\le
p\le1$ and $n^2p'\rightarrow\infty$ as $n\rightarrow\infty$. If
$G(n,n,p')$ a.a.s. has global resilience at least
$(1/2-\varepsilon/4)n^2p'$ with respect to a monotone increasing
graph property, then $G(n,n,p)$ a.a.s. has global resilience at
least $(1/2-\varepsilon/2)n^2p$ with respect to the same property. \\
\normalfont This claim can be shown similarly as \cite[Proposition
3.1]{8} or \cite[Proposition 2.7]{1}. We omit the proof here.

Therefore, another application of Theorem 3 implies that if $p\ge
Cn^{-2/3}$, then a.a.s. every subgraph $G'\subset G(n,n,p)$ with
more than $(1/2+\varepsilon/2)n^2p$ edges satisfies
$\mathcal{P}_{\delta}$. Note that every Hamiltonian graph with
property $\mathcal{P}_{\delta}$ contains a cycle of length $t$ for
all even $t\in[8,2\delta n]\cup[2(1-\delta)n,2n]$. The proof of
Theorem 6 is completed. $\Box$.

\section{Existence of medium cycles of even lengths}

In this section we establish the following result, which together
with Theorem 6 readily gives our main result Theorem 2.

\smallskip
\noindent\textbf{Theorem 7.} \quad \itshape For any
$\varepsilon\in(0,1)$ and $\delta\in(0,1)$, there exists a constant
$C>0$ such that if $p\ge Cn^{-2/3}$, then $G(n,n,p)$ a.a.s.
satisfies the following. Every Hamiltonian subgraph $G'\subset
G(n,n,p)$ with more than $(1+\varepsilon)n^2p/2$ edges contains a
cycle of length $t$ for all even $t\in[2\delta n,2(1-\delta)n]$.
\normalfont
\smallskip

Fix an allowable labeling of the vertices of $K_{n,n}$ with $[2n]$.
We partition the edge set $E(K_{n,n})$ as
$$
E(K_{n,n})=\cup_{i=0}^{n-1}E_i,
$$
where $E_i=\{\{x,y\}:x+y\equiv 2i+1\mod 2n\}$. For each $0\le i\le
n-1$, we define an ordering $\le_i$ for the elements of $E_i$ as
follows. If we evenly place the numbers from $[2n]$ on a circle
(i.e., $C_{2n}$ as defined above), each set $E_i$ will comprise all
parallel edges in some direction. We order them as per their
distance from the minor arc connecting $i$ to $i+1$; c.f. Fig. 3. We
refer to the elements of $E_i$ as the edges in direction $i$. Note
that $|E_i|=n$ for all $i$. Two edges $e_1,e_2\in
E(K_{n,n})\backslash E(C_{2n})$ is said to be crossing if their
endpoints are all distinct and lie alternately on the cycle
$C_{2n}$.

\begin{figure}[!t]
\begin{center}
\scalebox{0.7}{\includegraphics[85pt,353pt][559pt,561pt]{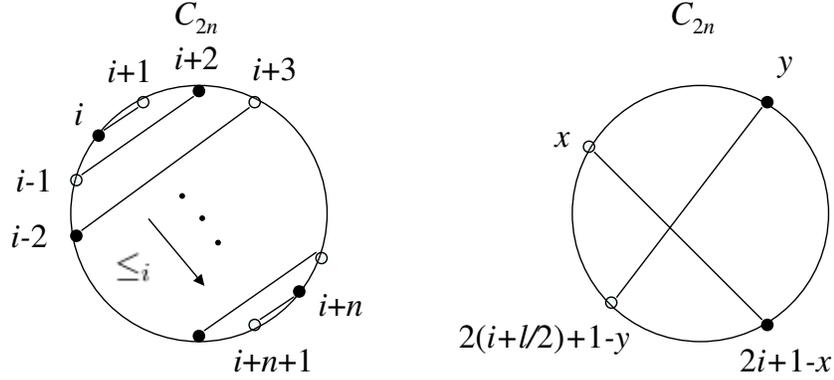}}\caption{The
total order set $E_i$ with arrow pointing from $\le_i$-smaller to
$\le_i$-larger elements; two crossing edges $e_1=\{x,2i+1-x\}$ and
$e_2=\{y,2(i+l/2)+1-y\}$.}
\end{center}
\end{figure}

In what follows, we will still adopt a similar reasoning as
conducted in \cite{1}. We want to highlight a remarkable
modification in our methodology.
\begin{itemize}
\item The redefinition of total order sets $E_i$ $(i=0,1,\cdots,n-1)$ allows
a smooth switch to the bipartite structure. This partition appears
to be critical in the following development.
\end{itemize}

We observe that for every $i\in[2n]$ and every even $l$ with $2\le
l\le 2n-2$, the graph $C_{2n}\cup\{e_1,e_2\}$, where $e_1\in E_i$
and $e_2\in E_{i+l/2}$ are crossing edges, contains cycles of
lengths $l+2$ and $2n-l+2$; see Fig. 3.

For $\beta\in(0,1/6)$, and $k\in\mathbb{N}$ with $1\le k\le
2(1/2-\beta)n$, we define $E_i^k$ be the set of $2\beta n$
consecutive (with respect to $\le_i$) edges in $E_i$, beginning from
the $k$-th smallest element of $E_i$. We may refer to $E_i^k$ as an
interval of length $2\beta n$, whose leftmost endpoint is the $k$-th
smallest element of $E_i$. Denote by $M_i\subset E_i$ the set of
$2(1/2-2\beta)n$ middle elements of $E_i$ without the leftmost and
rightmost intervals of lengths $2\beta n$. Let $G$ be a $2n$-vertex
bipartite balanced graph with an allowable labeling with $[2n]$. For
$i\in[2n]$, $\varepsilon'\in(0,1)$ and $p\in[0,1]$, we say that the
direction $E_i$ is $(\beta,\varepsilon',p)$-good in $G$ if for all
$k\in\mathbb{N}$ with $1\le k\le 2(1/2-\beta)n$, $G$ satisfies
\begin{equation}
||E(G)\cap E_i^k|-2\beta np|\le 2\varepsilon'\beta np,\label{5}
\end{equation}
and
\begin{equation}
||E(G)\cap M_i|-2(1/2-2\beta)np|\le
2\varepsilon'(1/2-2\beta)np.\label{6}
\end{equation}

\smallskip
\noindent\textbf{Lemma 4.} \quad \itshape Let
$\beta,\varepsilon'\in(0,1/6)$. If $p\ge Cn^{-2/3}$ for some $C>0$,
then a.a.s. for every allowable labeling of vertices of $G(n,n,p)$
with $[2n]$, there are at most $n^{5/6}$ directions that are not
$(\beta,\varepsilon',p)$-good in $G(n,n,p)$. \normalfont
\smallskip

\noindent\textbf{Proof.} Let $G$ be a graph drawn from $G(n,n,p)$
and fix an allowable labeling of the vertices of $G$ with $[2n]$. It
follows from Theorem 3 that, for all $i$ and $k$,
$$
P\left(||E(G)\cap E_i^k|-2\beta np|>2\varepsilon'\beta
np\right)\le2e^{-2\varepsilon'^2\beta np/3}\le e^{-cnp}
$$
and
$$
P\left(||E(G)\cap M_i|-2(1/2-2\beta)np|>
2\varepsilon'(1/2-2\beta)np\right)\le2e^{-2\varepsilon'^2(1/2-2\beta)
np/3}\le e^{-cnp},
$$
where $c=c(\beta,\varepsilon')>0$. Therefore, given $i$, $P(E_i\
\mathrm{is}\ \mathrm{not}\
(\beta,\varepsilon',p)\mbox{-}\mathrm{good})\le4(1/2-\beta)ne^{-cnp}\le
e^{-cnp/2}$ by using (\ref{5}) and (\ref{6}). Since the events
$\{E_i\ \mathrm{is}\ \mathrm{not}\
(\beta,\varepsilon',p)\mbox{-}\mathrm{good}\}_{0\le i\le n-1}$ are
mutually independent, the probability that there are more than
$n^{5/6}$ not good directions is at most
$$
{n \choose
n^{5/6}}\left(e^{-cnp/2}\right)^{n^{5/6}}\le2^ne^{-cn^{11/6}p/2}\le
e^{-c'n^{7/6}},
$$
where $c'=c'(c,C)>0$. Since there are $(n!)^2$ different allowable
labelings, the probability of there being an allowable labeling with
more than $n^{5/6}$ not good directions is at most
$$
(n!)^2e^{-c'n^{7/6}}\le\frac{e^2n^{2n+1}}{e^{2n+c'n^{7/6}}}=o(1),
$$
as $n\rightarrow\infty$. The proof is completed. $\Box$

For $\beta\in(0,1/6)$, a crossing between two edges $\{x_1,y_1\}$
and $\{x_2,y_2\}$ is said to be close if
$$
\min\{\|x_1-x_2\|,\|x_1-y_2\|,\|y_1-x_2\|,\|y_1-y_2\|\}\le2\beta n.
$$
The following statements can be proved based on a similar
observation in \cite[Lemma 3.10]{1}. We leave the proof to the
reader.

\smallskip
\noindent\textbf{Lemma 5.} \quad \itshape For $i\in[2n]$,
$\beta\in(0,1/6)$ and even $l$ with $4\beta n+1\le l\le
(2-4\beta)n-1$, the following statements are true.
\begin{itemize}
\item[(i)] Every edge in $E_i$ forms close crossings with at most $4\beta
n$ edges from $E_{i+l/2}$, and these edges can be covered by a set
of the form $E_{i+l/2}^{k_1}\cup E_{i+l/2}^{k_2}$ for some $1\le
k_1,k_2\le2(1/2-\beta)n$.

\item[(ii)] At least $(1-4\beta)n$ edges in $E_i$ form close
crossings with exactly $4\beta n$ edges from $E_{i+l/2}$, and these
$4\beta n$ edges constitute a set of the form $E_{i+l/2}^{k_1}\cup
E_{i+l/2}^{k_2}$ for some $1\le k_1,k_2\le2(1/2-\beta)n$.

\item[(iii)] The $(1-4\beta)n$ edges in (ii) cover $M_i$.

\end{itemize}
\normalfont
\smallskip

\noindent\textbf{Completion of the proof of Theorem 7.} Let
$\varepsilon'=\varepsilon/17$ and
$\beta=\min\{\delta/3,\varepsilon'\}$. Let $G$ be a graph drawn from
$G(n,n,p)$. By virtue of Lemma 4 a.a.s. every allowable labeling of
the vertices of $G$ with $[2n]$ yields at most $\varepsilon'n$
directions that are not $(\beta,\varepsilon',p)$-good in $G$. It
follows from Theorem 3 that $e(G)\le(1+\varepsilon/4)n^2p$ a.a.s.
Fix an even $t\in[2\delta n,2(1-\delta)n]$. It suffices to show
that, conditioned on the above two events, every Hamiltonian
subgraph $G'\subset G$ with more than $(1+\varepsilon)n^2p/2$ edges
contains a $t$-cycle.

Fix such a subgraph $G'$ and an allowable labeling of the vertices
such that $C_{2n}$ is a Hamilton cycle in $G'$, and set $l=t-2$.
Based on our above observation, we only need to show that for some
$i\in[2n]$, the graph $G'$ contains two edges $e_1\in E_i$ and
$e_2\in E_{i+l/2}$ which form a close crossing.

Denote by $I$ the set of directions that are
$(\beta,\varepsilon',p)$-good in $G$. Hence,
$|I|\ge(1-\varepsilon')n$ by our condition. Let $X$ be the number of
close crossings between pairs of edges in $G$ which came from $E_i$
and $E_{i+l/2}$ satisfying $i,i+l/2\in I$. In the following, we
assume $l\not=n$ (if $n$ is odd, it clearly holds; if $n$ is even,
the proof is similar and we leave it to the reader). Under this
assumption, we have $E_{i+l/2}\not=E_{i-l/2}$. So, the number of
pairs $\{i,i+l/2\}\subset I$ is at least $(1-2\varepsilon')n$. Fix
any such pair of them. By the definition of $\beta$ and $t$, we
obtain that $l\in(5\beta n, (2-5\beta)n)$. From Lemma 5 (ii), (iii)
and (\ref{6}) we know that each of the at least
$2(1-\varepsilon')(1/2-2\beta)np$ edges in $M_i\cap E(G)$ forms a
close crossing with every edge from some two disjoint sets
$E_{i+l/2}^k$ of size $2\beta n$ each. Recall that $i+l/2\in I$. By
(\ref{5}), the graph $G$ contains at least $2(1-\varepsilon')\beta
np$ edges in each such set $E_{i+l/2}^k$. Consequently, we obtain
\begin{eqnarray}
X&\ge&(1-2\varepsilon')n\cdot
2(1-\varepsilon')\left(\frac12-2\beta\right)np\cdot2\cdot2(1-\varepsilon')\beta
np\nonumber\\
&\ge&4(1-4\varepsilon'-4\beta)\beta n^3p^2.\label{7}
\end{eqnarray}
It follows from Lemma 5 (i) that each edge $e_1\in E_i$ forms close
crossings with at most $4\beta n$ edges from $E_{i\pm l/2}$, and
these edges are covered by some sets $E_{i\pm l/2}^k$. Therefore, by
using (\ref{5}), every edge in a $(\beta,\varepsilon',p)$-good
direction $i$ forms at most $8(1+\varepsilon')\beta np$ close
crossings with edges in a $(\beta,\varepsilon',p)$-good direction
$i\pm i/2$. Let $Y$ be the number of crossings in $G$ that are
counted by $X$ but not contained in $G'$. We have
\begin{eqnarray}
Y&\le&(e(G)-e(G'))\cdot8(1+\varepsilon')\beta np\nonumber\\
&\le&\left(\left(1+\frac{\varepsilon}{4}\right)n^2p-(1+\varepsilon)\frac{n^2p}{2}\right)\cdot8(1+\varepsilon')\beta np\nonumber\\
&\le&(4-2\varepsilon+\varepsilon')\beta n^3p^2.\label{8}
\end{eqnarray}
By our definitions, we have $16\beta+17\varepsilon'<2\varepsilon$.
Hence, we derive $X>Y$ by (\ref{7}) and (\ref{8}). In other words,
$G'$ contains two edges from $E_i$ and $E_{i+l/2}$ that form a close
crossing. This finally completes the proof of Theorem 7. $\Box$


\end{document}